# ON TESTING THE SIGNIFICANCE OF SETS OF GENES

By Bradley Efron[1] and Robert Tibshirani[2]

*Stanford University*

This paper discusses the problem of identifying differentially expressed groups of genes from a microarray experiment. The groups of genes are externally defined, for example, sets of gene pathways derived from biological databases. Our starting point is the interesting Gene Set Enrichment Analysis (GSEA) procedure of Subramanian et al. [*Proc. Natl. Acad. Sci. USA* **102** (2005) 15545–15550]. We study the problem in some generality and propose two potential improvements to GSEA: the *maxmean* statistic for summarizing gene-sets, and *restandardization* for more accurate inferences. We discuss a variety of examples and extensions, including the use of gene-set scores for class predictions. We also describe a new R language package *GSA* that implements our ideas.

**1. Introduction.** We discuss the problem of identifying differentially expressed groups of genes from a set of microarray experiments. In the usual situation we have $N$ genes measured on $n$ microarrays, under two different experimental conditions, such as control and treatment. The number of genes $N$ is usually large, say, at least a few thousand, while the number samples $n$ is smaller, say, a hundred or fewer. This problem is an example of multiple hypothesis testing with a large number of tests, one that often arises in genomic and proteomic applications, and also in signal processing. We focus mostly on the gene expression problem, but our proposed methods are more widely applicable.

Most approaches start by computing a two-sample $t$-statistic $z_j$ for each gene. Genes having $t$-statistics larger than a pre-defined cutoff (in absolute value) are declared significant, and then the family-wise error rate or false discovery rate of the resulting gene list is assessed by comparing the tail area

Received October 2006; revised January 2007.

[1]Supported in part by NSF Grant DMS-05-05673 and National Institutes of Health Contract 8RO1 EB002784.

[2]Supported in part by NSF Grant DMS-99-71405 and National Institutes of Health Contract N01-HV-28183.

from a null distribution of the statistic. This null distribution is derived from data permutations, or from asymptotic theory.

In an interesting and useful paper, Subramanian et al. (2005) proposed a method called *Gene Set Enrichment Analysis* (GSEA) for assessing the significance of pre-defined gene-sets, rather than individual genes. The gene-sets can be derived from different sources, for example, the sets of genes representing biological pathways in the cell, or sets of genes whose DNA sequences are close to together on the cell's chromosomes. The idea is that these gene-sets are closely related and, hence, will have similar expression patterns. By borrowing strength across the gene-set, there is potential for increased statistical power. In addition, in comparing study results on the same disease from different labs, one might get more reproducibility from gene-sets than from individual genes, because of biological and technical variability.

The GSEA methods works roughly as follows. We begin with a pre-defined collection of gene-sets $\mathcal{S}_1, \mathcal{S}_2, \ldots, \mathcal{S}_K$. We compute $t$-statistic $z_j$ for all $N$ genes in our data. Let $\mathbf{z}_k = (z_1, z_2, \ldots, z_m)$ be the gene scores for the $m$ genes in gene-set $\mathcal{S}_k$. In GSEA we then compute a gene-set score $S_k(\mathbf{z}_k)$ for each gene-set $\mathcal{S}_k$, equal to essentially a signed version of the Kolmogorov–Smirnov statistic between the values $\{z_j, j \in \mathcal{S}_k\}$ and their complement $\{z_j, j \notin \mathcal{S}_k\}$; the sign taken positive or negative depending on the direction of shift. The idea is that if some or all of the gene-set $\mathcal{S}_k$ have higher (or lower) values of $z_j$ than expected, their summary score $S_k$ should be large. An absolute cutoff value is defined, and values of $\mathcal{S}_k$ above (or below) the cutoff are declared significant. The GSEA method then does many permutations of the sample labels and recomputes the statistic on each permuted dataset. This information is used to estimate the false discovery rate of the list of significant gene-sets. The Bioconductor package `limma` offers an analysis option similar to GSEA, but uses instead the simple average of the scores $z_k$ [see Smyth (2004)], randomizing over the set of genes rather than the set of samples, what we call "row randomization" here.

Other related ideas may be found in Pavlidis et al. (2002) and Rahnenfhrer et al. (2004). Nobel and Wright (2005) propose the "SAFE" methodology a quite general permutation-based approach to the enrichment testing problem. Newton et al. (2006) propose random set scoring methods for assessing the significance of gene-set enrichment. Zahn et al. (2006) propose an alternative to GSEA that uses a Van der Waerden statistic in place of the Kolmogorov–Smirnov statistic and bootstrap sampling of the arrays instead of a permutation distribution. Other papers that address the problem of testing for differentially expressed sets of genes include Szabo et al. (2003), Frisina et al. (2004), Lu et al. (2005) and Dettling et al. (2005). One of our goals here is to make explicit the choices involved between the various randomization schemes.

In studying the GSEA work, we have found some shortcomings and ways it could be improved. The GSEA's dependence on Kolmogorov–Smirnov statistics is a reasonable choice, but not a necessary one. This paper puts the GSEA procedure in a more theoretical framework that allows us to investigate questions of efficiency for gene-set inference; a new procedure based on the "maxmean" statistic is suggested that has superior power characteristics versus familiar location/scale alternatives.

Here are two simulated data examples that illustrate some of the main issues, and allow us to introduce our proposed solution. We generated data on 1000 genes and 50 samples, with each consecutive nonoverlapping block of 20 genes considered to be a gene-set. The first 25 samples are the control group, and the second 25 samples are the treatment group. First we generated each data value as i.i.d. $N(0,1)$. Then the constant 2.5 was added to the first 10 genes in the treatment group. Thus, half of the first gene-set (first block of 20 genes) has a higher average expression in the treatment group, while all other gene-sets have no average difference in the two groups.

The left panel of Figure 1 shows a histogram (black lines) of the GSEA scores for the 50 gene-sets. The first gene-set clearly stands out, with a value of about 0.9. We did 200 permutations of the control-treatment labels, producing the dashed histogram in the top left panel of Figure 1. The first gene-set stands out on the right side of the histogram. So the GSEA method has performed reasonably well in this example.

In the paper we study alternative summary statistics for gene-sets. Our favorite is something we call the "maxmean statistic": we compute the average of the positive parts of each $z_i$ in $\mathcal{S}$, and also the negative parts, and choose the one that is larger in absolute value. The results for maxmean in this example are shown in the right panel of Figure 1. The first gene-set stands out more clearly than it does in the left panel. In this paper we show by both analytic calculations and simulations that the maxmean statistics are generally more powerful than GSEA.

Now consider a different problem. We generated data exactly as before, except that the first 10 genes in *every* gene-set are 2.5 units higher in the treatment group. The top left panel of Figure 2 shows a histogram of the maxmean scores for the 50 gene-sets, and a histogram of the scores from 200 permutations of the sample labels (dashed). All of the scores look significantly large compared to the permutation values. But given the way that the data were generated, there seems to be nothing special about any one gene-set. To quantify this, we "row randomized" the 1000 genes, leaving the sample labels as is. The first 20 genes in the scrambled set became the first gene-set, the second 20 genes became the second gene-set, and so on. We did this many 200 times, recomputing the maxmean statistic on each scrambled set. The results are shown in the bottom left panel of Figure 2. None of the original 50 gene-sets are notable when compared to gene-sets formed by randomly sampling from the full set of genes.

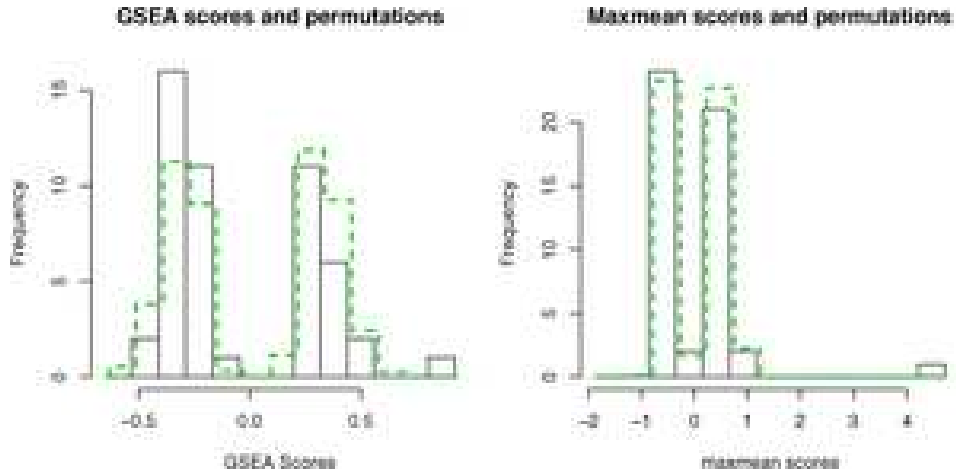

Fig. 1. *Example* 1: *left panel shows histograms of GSEA scores for original data (black) and from 200 permutations of the sample labels (dashed). The first gene-set stands our fairly clearly with a score of about 0.7. The right panel shows the results using instead the "maxmean" statistic. The first gene-set stands out much more clearly than in the left panel.*

The point is that any method for assessing gene-sets should compare a given gene-set score not only to scores from permutations of the sample labels, but also take into account scores from sets formed by random selections of genes. Both GSEA and SAFE do something like this, at least implicitly.

The bottom right panel of Figure 2 puts all of these ideas together. It uses the *restandardized* version of the maxmean statistic, in which we center and scale the maxmean statistic by its mean and standard deviation under the row randomizations, like those in the bottom left panel. This standardized maxmean statistic is then computed both on the original data (light histogram) and on each of the permuted datasets (dark histogram). As we would expect, we see no significant gene-sets. This restandardization is potentially important for any gene summary statistic: it turns out that the GSEA statistic incorporates a form of what we are calling restandardization.

The two ideas illustrated above—alternative summary statistics for gene-sets and restandardization based on row randomization—are two of the main proposals in this paper. In Section 2 we describe the randomization and permutation methods for estimating appropriate null distributions. Section 3 studies the choice of summary statistic for gene-sets, and introduces the *maxmean* statistic. In Section 4 we summarize our proposal for Gene Set Analysis and discuss computational issues. A simulation study is carried out in Section 5, comparing the power of the maxmean statistic to the GSEA statistic and other competitors. Finally, in Section 6 we give a number of examples of the method, including applications to a kidney cancer dataset, a generalization to the class prediction problem and comparison of different gene-set collections over the same expression dataset.

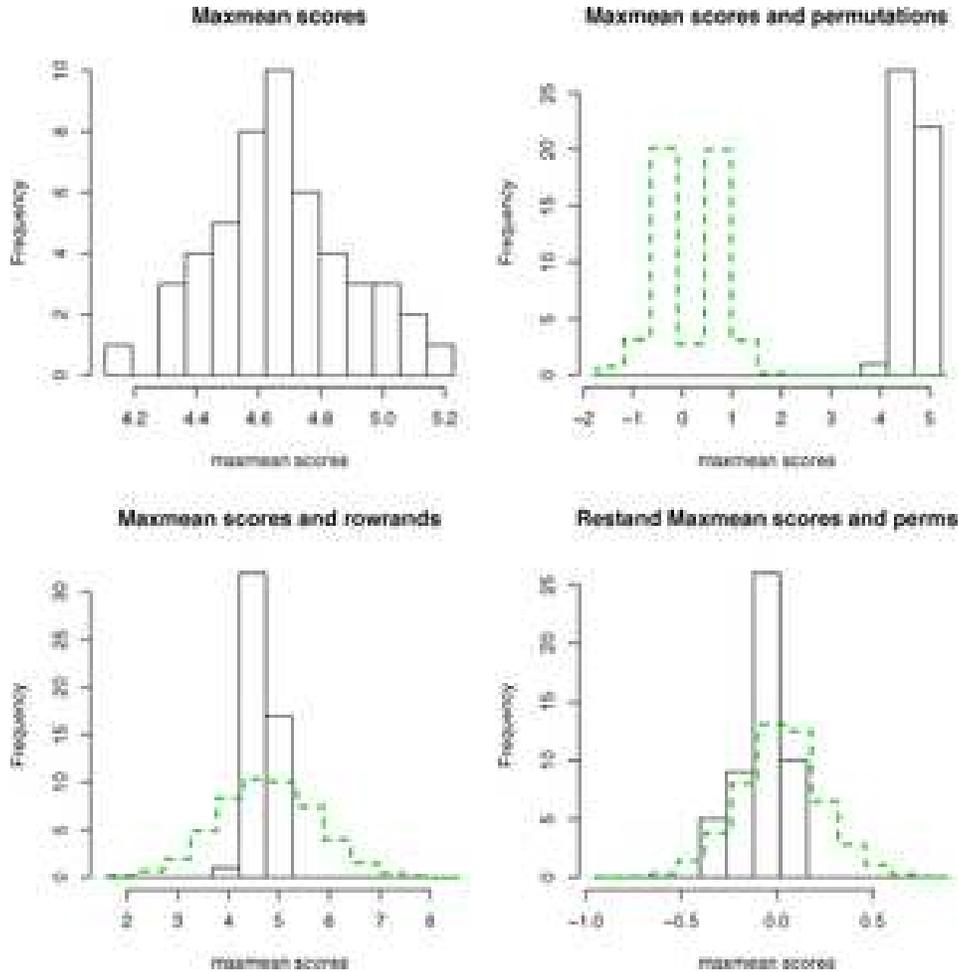

FIG. 2. *Example* 2: *top left panel shows a histogram of the maxmean scores. In the top right panel a histogram of maxmean scores from 200 permutations of the sample labels has been added (dashed): all of the gene-sets look to be significant. In the bottom left panel a histogram of the maxmean scores is shown along with maxmean scores from 200 "row randomizations" —gene-sets chosen from the full collection of genes at random (dashed). Finally, in the bottom right panel, histograms of the restandardized scores from the original data (solid) and from 200 permutations of the sample labels (dashed) are shown. (Note that the horizontal axis has changed, since the values have been centered and scaled). None of the gene-sets looks significant, which is reasonable, given the way the data were generated.*

**2. Randomization and permutation.** A straightforward approach to gene-set analysis begins with some enrichment score $S$ and computes its significance by comparison with permutation values $S^*$. Here we argue that a second kind of comparison operation, "row randomization," is also needed to avoid bias in the determination of significance. We begin with a simplified statement of the gene-set problem, leading later to a more realistic analysis.

Let $X$ indicate an $N$ by $n$ matrix of expression values, $N$ genes and $n$ microarrays, with the first $n_1$ columns of $X$ representing Class 1 and the last $n_2$ Class 2, $n_1 + n_2 = n$. In the $p53$ example of Subramanian et al. (2005) there are $N = 10100$ genes and $n = 50$ arrays relating to cell lines with normal or mutated states for the p53 factor, $n_1 = 17$ normal and $n_2 = 33$ mutated.

The $i$th row of $X$, that is, the data for gene $i$, yields a two-sample $t$-statistic "$t_i$" comparing the classes. In this section it will be convenient to transform the $t_i$ values to $z$-values "$z_i$," theoretically having a standard normal distribution,

$$\text{(2.1)} \qquad \text{theoretical null: } z_i \sim N(0,1),$$

under the null hypothesis of no difference between the two treatments. [The transformation is $z_i = \Phi^{-1}(F_{n-2}(t_i))$, where $\Phi$ is the standard normal cumulative distribution function (c.d.f.) while $F_{n-2}$ is the c.d.f. for a $t$ distribution having $n-2$ degrees of freedom.] The methods here apply outside the realm of genes, microarrays and $t$-tests, but that is the main application we have in mind.

To begin with, suppose that a single gene-set "$\mathcal{S}$," comprising $m$ genes, is under consideration, and we wish to test the hypothesis that $\mathcal{S}$ is *enriched*, meaning that its $m$ $z$-values are unusually large, positively or negatively, in a sense to be defined. The scientific idea underlying enrichment, as nicely stated in Subramanian et al. (2005), is that a biologically related set of genes, perhaps representing a genetic pathway, could reveal its importance through a general effect on its constituent $z$-values, whether or not the individual $z_i$'s were significantly nonzero.

Let $\mathbf{z}_\mathcal{S}$ indicate the set of $m$ $z$-values in $\mathcal{S}$ and

$$\text{(2.2)} \qquad S = S(\mathbf{z}_\mathcal{S})$$

define an enrichment test statistic, with larger value of $S$ indicating greater enrichment. The GSEA algorithm uses a Kolmogorov–Smirnov type function for $S(\mathbf{z}_\mathcal{S})$. A simpler approach starts with a function "$s$" of the individual $z$-values,

$$\text{(2.3)} \qquad s_i = s(z_i),$$

and takes the gene-set score $S$ to be the average of $s_i$ in $\mathcal{S}$,

$$\text{(2.4)} \qquad S = \sum_\mathcal{S} s(z_i)/m.$$

The choice $s(z) = |z|$ will be discussed later. *Limma* [Smyth (2004)], a microarray analysis program available in the Bioconductor R package, implements (2.4) with $s(z) = z$, so $S$ is simply $\bar{z}_\mathcal{S}$, the average $z$-score in $\mathcal{S}$ (actually using $S = |\bar{z}_\mathcal{S}|$ for two-sided inference). Section 3 develops a third choice, the "maxmean" statistic.

Testing $\mathcal{S}$ for enrichment requires a null hypothesis distribution for $S$, and that is where difficulties begin; there are two quite different models for what "null" might mean. We discuss these next, in a spirit closely related to "Q1" and "Q2" of Tian et al. (2005).

**Randomization Model.** The null hypothesis $H_{\text{rand}}$ is that $\mathcal{S}$ has been chosen by random selection of $m$ genes from the full set of $N$ genes. In this case the null density of $S$, say, $g_{\text{rand}}(S)$, can be obtained by *row randomization*: sets $\mathcal{S}^\dagger$ of $m$ rows of the expression matrix $X$ are drawn at random, giving randomized versions $S^\dagger$ of (2.4), a large number of which are generated and used to estimate $g_{\text{rand}}(S)$. Equivalently, random subsets $\mathbf{z}_\mathcal{S}^\dagger$ of size $m$ are drawn from all $N$ $z_i$'s, giving $S^\dagger = S(\mathbf{z}_\mathcal{S}^\dagger)$.

**Permutation Model.** Let $X_\mathcal{S}$ be the $m$ by $n$ submatrix of $X$ corresponding to $\mathcal{S}$. The null hypothesis $H_{\text{perm}}$ is that the $n$ columns of $X_\mathcal{S}$ are independent and identically distributed $m$-vectors (i.i.d.). The null density of $S$, $g_{\text{perm}}(S)$, is obtained by column permutations, leading to an estimate of $g_{\text{perm}}(S)$ in the usual way.

The Randomization Model has the appealing feature of operating *conditionally* on the set $\mathbf{z}$ of all $N$ $z$-values: given $\mathbf{z}$, it tests the null hypothesis that the observed $S$ is no bigger than what might be obtained by a random selection of $\mathcal{S}$.

Its defect concerns correlation between the genes. Suppose the scores $s_i$ in (2.3) have empirical mean and standard deviation

$$(2.5) \qquad s \sim (\text{mean}_s, \text{stdev}_s)$$

for the ensemble of all $N$ genes. Because row randomization destroys gene-wide correlations, the randomized version of (2.4) will have mean and standard deviation

$$(2.6) \qquad S^\dagger \sim (\mu^\dagger, \sigma^\dagger),$$

with

$$(2.7) \qquad \mu^\dagger = \text{mean}_s \quad \text{and} \quad \sigma^\dagger = \text{stdev}_s / \sqrt{m}.$$

However, $\sigma^\dagger$ will underestimate the variability of $S$ if there is positive correlation among the $z$-values in $\mathcal{S}$, which seems likely if $\mathcal{S}$ was chosen on the basis of biological similarities. The $p$-value in the Limma package compares $S = \bar{z}_\mathcal{S}$ with the row randomization distribution of $S^\dagger$, ignoring, perhaps dangerously, correlation effects.

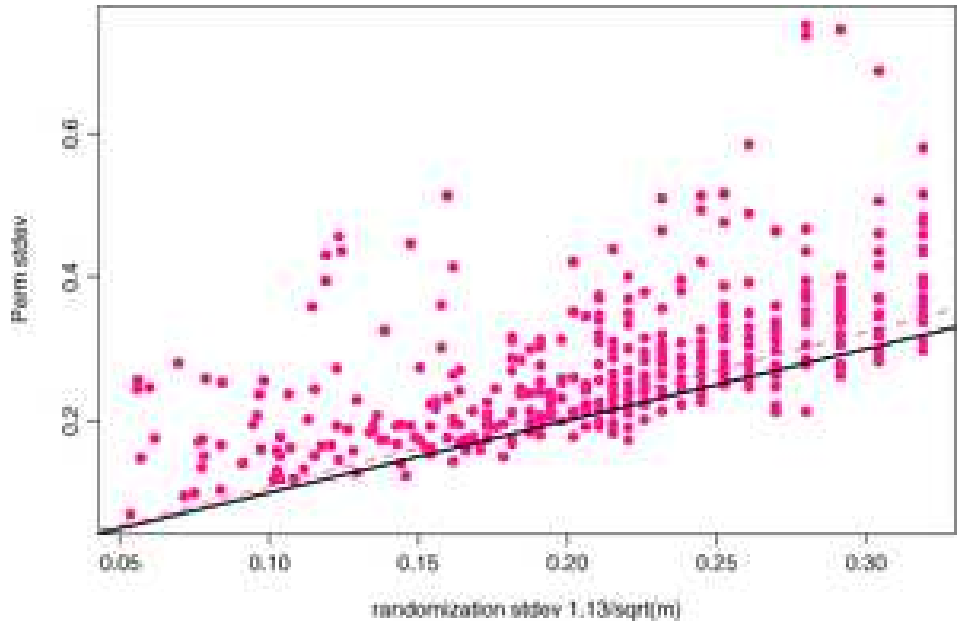

Fig. 3. *Permutation standard deviation compared to randomization standard deviation for $S = \bar{z}_\mathcal{S}$, 395 gene-sets for p53 data, Subramanian et al. (2005). Permutation values exceed randomization values, especially for gene-sets with large m at left, reflecting positive correlations within gene-sets.*

Permutation methods deal nicely with the correlation problem. Let $R_\mathcal{S}$ indicate the correlation matrix of the i.i.d. columns of $X_\mathcal{S}$ in the Permutation Model. Because column-wise permutations maintain the integrity of column vectors, it turns out, under some regularity conditions, that $\mathbf{z}^*_\mathcal{S}$, the permutation vector of $z$-values in $\mathcal{S}$, will have correlation matrix approximating $R_\mathcal{S}$, as will $\mathbf{z}_\mathcal{S}$ itself. In other words, for a prechosen gene-set $\mathcal{S}$, the permutation correlation matrix of $\mathbf{z}_\mathcal{S}$ will be approximately correct.

Figure 3 concerns correlation effects in the p53 example of Subramanian et al. (2005); the authors provide a catalog of 522 gene sets, of which 395 have $m \geq 10$ members. Taking $s(z) = z$ yields

$$(2.8) \qquad (\text{mean}_s, \text{stdev}_s) = (0.00, 1.13),$$

as in (2.5), this being the mean and standard deviation of all those $z_i$'s involved in the 395 gene-sets. 400 permutation values of $S^* = \bar{z}^*_\mathcal{S}$, the average $z^*_i$ values in $\mathcal{S}$, were generated for each of the 395 gene-sets $\mathcal{S}$, and the permutation standard deviation $\sigma^*_\mathcal{S}$ computed.

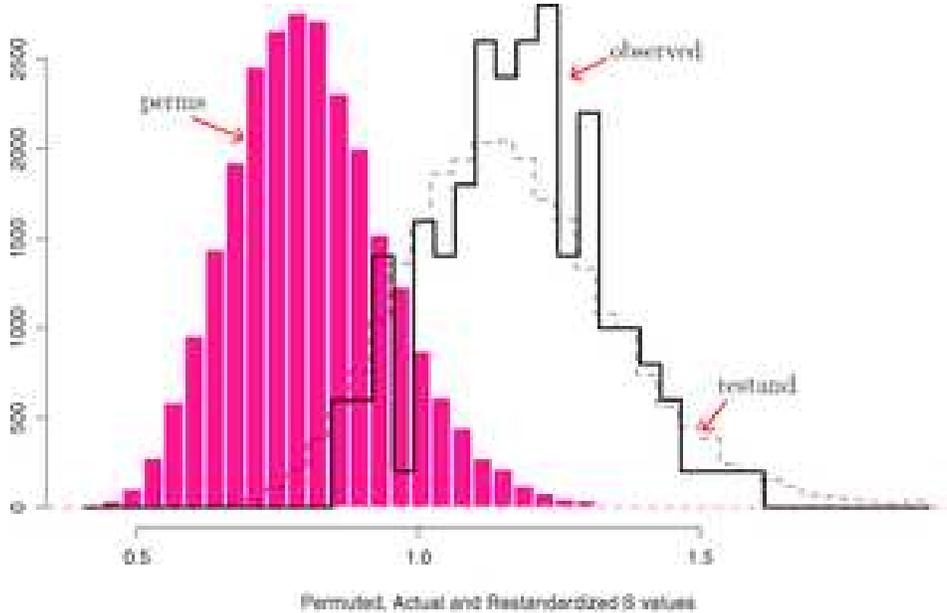

FIG. 4. *Enrichment statistics $S$ for 129 randomly selected gene-sets of size $m = 25$, heavy line histogram, lies far to right of permutation $S^*$ values, solid histogram. $S = \sum_{\mathcal{S}} |z_i|/m$, BRCA data, Hedenfalk et al. (2001). Restandardized values (2.12), light line histogram, gives a much better match to actual values.*

Figure 3 plots $\sigma_{\mathcal{S}}^*$ versus $\sigma_{\mathcal{S}}^{\dagger} = 1.13/\sqrt{m_{\mathcal{S}}}$, the randomization standard deviation (2.7); $\sigma_{\mathcal{S}}^*/\sigma_{\mathcal{S}}^{\dagger}$ has median 1.08 (light line), but the ratio is considerably greater for large $m$, with median 1.50 for $m \geq 50$.

The Permutation Model has a serious weakness of its own: it does not take into account the overall statistics $(\text{mean}_s, \text{stdev}_s)$, (2.5), as demonstrated next. Figure 4 concerns the "BRCA data," Efron (2004), a microarray experiment comparing two classes of breast cancer patients, $N = 3226$ genes, $n = 15$ microarrays, $n = 7$ for class 1, $n_2 = 8$ for class 2. A two-sample $t$-statistic $t_i$ was computed for each gene and converted to $z$-value, $z_i = \Phi^{-1}(F_{13}(t_i))$. The $N$ $z$-values have mean and standard deviation

$$z \sim (-0.026, 1.43). \tag{2.9}$$

In this case (2.1), the theoretical $N(0,1)$ null, substantially underestimates the $z$-values' variability, their histogram looking like an overdispersed but somewhat short-tailed normal.

For this example 129 gene-sets $\mathcal{S}$ of size $m = 25$ each were constructed by random selection from the $N$ genes. Enrichment was tested with $S = \sum_{\mathcal{S}} |z_i|/m$, that is, by using (2.3) and (2.4) with $s(z) = |z|$. Figure 4 shows that the permutation distribution $S^*$ greatly underestimates the actual 129 $S$ values. A standard Benjamini–Hochberg (1995) FDR analysis based on permutation $p$-values, carried out at rejection level $q = 0.10$, labels 113 of the 129 generate as "enriched," even though we know that all of them were constructed randomly.

It is easy to spot the trouble here: the permutation $z$-values are much less dispersed than the actual $z$'s, with marginal mean and standard deviation, for all $N$ of them,

$$(2.10) \qquad z^* \sim (0.022, 1.014)$$

compared to (2.9). [This is a dependable feature of the permutation $t$-test, which tends to yield $z_i^* \stackrel{.}{\sim} N(0,1)$, the theoretical null, under a wide variety of circumstances.]

The individual scores $s_i^* = s(z_i^*)$ are correspondingly reduced. Letting (mean$^*$, stdev$^*$) indicate the mean and standard deviation of $s^*$ over all $N$ genes and a large number of permutations,

$$(2.11) \qquad s^* \sim (\text{mean}^*, \text{stdev}^*),$$

we computed (mean$^*$, stdev$^*$) = $(0.82, 0.60)$ compared to (mean$_s$, stdev$_s$) = $(1.17, 0.82)$ for the actual scores (2.5). This translates into $S^*$ statistics that are much smaller than the actual $S$ values.

2.1. *Restandardization.* We combine the randomization and permutation ideas into a method that we call "**Restandardization**." This can be thought of as a method of correcting permutation values $S^*$ to take account of the overall distribution of the individual scores $s_i$. For a gene-set score $S = \sum_{\mathcal{S}} s(z_i)/m$ as in (2.3) and (2.4), the *restandardized* permutation value is defined to be

$$(2.12) \qquad S^{**} = \text{mean}_s + \frac{\text{stdev}_s}{\text{stdev}^*}(S^* - \text{mean}^*),$$

with (mean$_s$, stdev$_s$) and (mean$^*$, stdev$^*$) the overall means and standard deviations defined at (2.5) and (2.11).

A restandardized $p$-value for testing the enrichment of gene-set $\mathcal{S}$ is obtained by comparing its actual enrichment score $S$ to a large number "$B$" of $S^{**}$ simulations,

$$(2.13) \qquad p_{\mathcal{S}} = \#\{S^{**} \text{ values exceeding } S\}/B.$$

Figure 4 shows the restandardized value $S^{**}$ giving much better agreement than $S^*$ with the actual $S$ values, though the $S^{**}$ distribution is moderately overdispersed.

Here is another equivalent way to express the restandardization procedure: we first standardize the statistic $S$ with respect to its randomization mean and standard deviation, and then compute the permutation distribution of the standardized statistic. This gives the $p$-value

$$(2.14) \qquad p_{\mathcal{S}} = \#\left\{\frac{S^* - \text{mean}^*}{\text{stdev}^*} \text{ exceeding } \frac{S - \text{mean}_s}{\text{stdev}_s}\right\}\bigg/B.$$

It is this form of restandardization that we use in our software package GSA, described in Section 4.

Restandardization can also be applied to more complicated gene-set statistics $S$, not simple averages as in (2.4), via

$$(2.15) \qquad S^{**} = \mu^\dagger + \frac{\sigma^\dagger}{\sigma^*}(S^* - \mu^*).$$

Here $(\mu^\dagger, \sigma^\dagger)$ are the mean and standard deviation of $S^\dagger$ for a randomly selected gene-set of size $m$, as in (2.5), while $(\mu^*, \sigma^*)$ are the corresponding quantities based on a permuted data matrix; computing $(\mu^*, \sigma^*)$ requires a nested simulation.

The GSEA enrichment test incorporates a form of restandardization: it compares the empirical c.d.f. of the $z$-values in $\mathcal{S}$ with the c.d.f. of all other $z$-values; the latter c.d.f. brings to bear information from the overall distribution of $\mathbf{z}$, much like $\text{mean}_s$ and $\text{stdev}_s$ do in (2.14). This same approach could be used in context (2.4) by replacing $s(z_i)$ in (2.3) with $t(z_i) = (s(z_i) - \text{mean}_s)/\text{stdev}_s$, leading to

$$(2.16) \qquad T = \sum_{\mathcal{S}} t(z_i)/m = (S - \text{mean}_s)/\text{stdev}_s$$

as the enrichment statistic for $\mathcal{S}$. Then (2.14) reduces to

$$(2.17) \qquad p_{\mathcal{S}} = \#\{T^* > T\}/B,$$

so that the restandardized $p$-value equals the usual permutation $p$-value.

Restandardization can be shown to yield improved inferences in a variety of situations:

- If $\mathcal{S}$ was selected randomly, as in the Randomization Model.
- If the theoretical null (2.1), $z \sim N(0,1)$, agrees with the empirical distribution of the $N$ $z$-values.
- If the $z_i$'s are uncorrelated.

The method is not perfect, as examples can be constructed to show. Neither the Randomization nor Permutation Models perfectly describe how gene-sets $\mathcal{S}$ come under consideration in the catalog examples of Subramanian et al. (2005), making some sort of compromise formula a necessity.

**3. Efficient tests for enrichment.** The Randomization and Permutation Models specify null hypotheses for the selection of a gene-set $\mathcal{S}$. Notions of efficient testing, however, require us to specify alternatives to the null selection. This section begins with an exponential family selection model, considers some specific test statistics, and goes on to propose the "maxmean" statistic, which leads to enrichment tests with good power characteristics.

The *Poisson Selection Model* starts with independent Poisson variates

$$(3.1) \qquad I_i \stackrel{\text{ind}}{\sim} Po(\nu_i), \qquad \nu_i = \alpha e^{\beta' s_i}/T_\beta$$

for $i = 1, 2, \ldots, N$, where $s_i = s(z_i)$, for $s(\cdot)$ a given $J$-dimensional function, $\beta$ an unknown $J$-dimensional parameter vector,

$$(3.2) \qquad T_\beta = \sum_{i=1}^{N} e^{\beta' s_i},$$

and $\alpha = \sum_1^n \nu_i$ an unknown scalar parameters. In what follows, the $\nu_i$ will be small and the $I_i$ almost all zeros or ones, mostly zeros. For convenient exposition we assume that *all* the $I_i$ are zeros or ones, though this is not essential.

We define the selected gene-set $\mathcal{S}$ as composed of those genes having $I_i = 1$,

$$(3.3) \qquad \mathcal{S} = \{i : I_i = 1\},$$

so $\mathcal{S}$ has

$$(3.4) \qquad m \equiv \sum_1^N I_i$$

members. It is easy to calculate the probability $g_{\alpha,\beta}(\mathcal{S})$ of selecting my particular gene-set $\mathcal{S}$,

$$(3.5) \qquad g_{\alpha,\beta}(\mathcal{S}) = \left(\frac{e^{-\alpha}\alpha^m}{m!}\right)(m! e^{m[\beta' \bar{s}_\mathcal{S} - \log T_\beta]}),$$

where

$$(3.6) \qquad \bar{s}_\mathcal{S} = \sum_\mathcal{S} s_i/m.$$

This is a product of two exponential family likelihoods, yielding maximum likelihood estimates $\hat{\alpha} = m$ and $\hat{\beta}$ defined by

$$(3.7) \qquad \sum_{i=1}^{N} s_i e^{\hat{\beta}' s_i} \bigg/ \sum_{i=1}^{N} e^{\hat{\beta} s_i} = \bar{s}_\mathcal{S}.$$

Parameter $\beta = 0$ corresponds to the Randomization Model null hypothesis that $\mathcal{S}$ has been selected at random. Given gene-set size $m$, the sufficient statistic for testing $H_0: \beta = 0$ is $\bar{s}_\mathcal{S}$, called "S" at (2.4). An efficient

ON TESTING THE SIGNIFICANCE OF SETS OF GENES 13test rejects $H_0$ if $\bar{s}_\mathcal{S}$ is far from its null expectation $\bar{s} = \sum_1^N s_i/N$. If $s_\mathcal{S}$ is one-dimensional, as in Section 2, we reject for extreme values of $\bar{s}_\mathcal{S}$, either positive or negative.

Under $H_0: \beta = 0$, all gene-sets $\mathcal{S}$ of size $m$ are equally likely. Nonzero $\beta$ "tilts" the selection toward gene-sets having larger values of $\beta'\bar{s}_\mathcal{S}$, as seen in (3.5) where

$$(3.8) \qquad g_\beta(\mathcal{S}|m) = m! e^{m[\beta'\bar{s}_\mathcal{S} - \log T_\beta]}.$$

Efficient choice of a gene-set enrichment test statistic $S = \bar{z}_\mathcal{S}$ depends on the individual scoring function $s_i = s(z_i)$ in (3.3). Consider two choices:

$$(3.9) \qquad s^{(1)}(z) = z \quad \text{and} \quad s^{(2)}(z) = |z|,$$

both of which are univariate and have $\beta$ a scalar, $J = 1$. The first of these, the *Limma* choice, uses $S = \bar{z}_\mathcal{S}$; enriched gene-sets $\mathcal{S}$ will tend to have $\bar{z}_\mathcal{S}$ far away from $\bar{z} = \sum_1^N z_i/N$. This is fine if enrichment manifests itself as a *shift* of location in $z$-values. The second choice $s^{(2)}(z) = |z|$ has power against *scale* alternatives.

For general use one would like enrichment test statistics having power against both shift and scale alternatives. To this end, define the two-dimensional scoring function

$$(3.10) \qquad s(z) = (s^{(+)}(z), s^{(-)}(z)), \begin{cases} s^{(+)}(z) = \max(z, 0), \\ s^{(-)}(z) = -\min(z, 0), \end{cases}$$

and the *maxmean test statistic*

$$(3.11) \qquad S_{\max} = \max\{\bar{s}_\mathcal{S}^{(+)}, \bar{s}_\mathcal{S}^{(-)}\}.$$

$S_{\max}$ is designed to detect unusually large $z$-values in either or both directions. Note that $S_{\max}$ is not the same as the largest absolute value of the mean of the positive or negative $z$ values in the gene-set: it divides by the total number of genes $m$. If, for example, a gene-set of 100 genes had 99 scores of $-0.5$ and one score of 10, the average of the positive scores would be 10, and the average of the negative scores equal to $-0.5$, so the mean with largest absolute value would be 10. But the maxmean statistic is the maximum of $10/100 = 0.1$ and $-99(-0.5)/100 = 0.495$, and the negative scores would win. The maxmean statistic is robust by design, not allowing a few large positive or negative gene scores to dominate.

Figure 5 compares the power of $S_{\max}$ with the limma test that rejects the null hypothesis of no enrichment for large values of the absolute mean statistic $|\bar{z}_\mathcal{S}|$. Here $\mathcal{S}$ consists of $m = 25$ independent normal $z$-values,

$$(3.12) \qquad z_i \stackrel{\text{ind}}{\sim} N(b, g^2), \qquad i = 1, 2, \ldots, m = 25.$$

Contours of equal power are shown, for testing $H_0 : (b, g) = (0, 1)$ versus alternatives $(b \geq 0, g \geq 1)$ [or, equivalently, $(b \leq 0, g \geq 1)$], at level 0.95; better power is reflected by contours nearer the null hypothesis point $(b, g) = (0, 1)$.

As expected, the absolute mean test has good power for shift alternatives $b \neq 0$, but no power for scale alternatives $g > 1$, as shown by its nearly vertical contour lines. The maxmean test shows reasonable power in both directions. Comparing the 0.5 power contours, for instance, shows maxmean almost dominating $|\bar{z}_\mathcal{S}|$. Also shown is "KS," the 0.5 power contour for a test based on absolute Kolmogorov–Smirnov distance between the empirical c.d.f. of the 25 $z$-values and a $N(0,1)$ c.d.f. (The GSEA test is based on an improved variant of the Kolmogorov–Smirnov statistic.) In this case, maxmean dominates $ks$, its contour being everywhere closer to $(b, g) = (0, 1)$. Finally, we show the results for the test based on the mean absolute value: this has the best power for pure scale alternatives but is not very sensitive to location shifts.

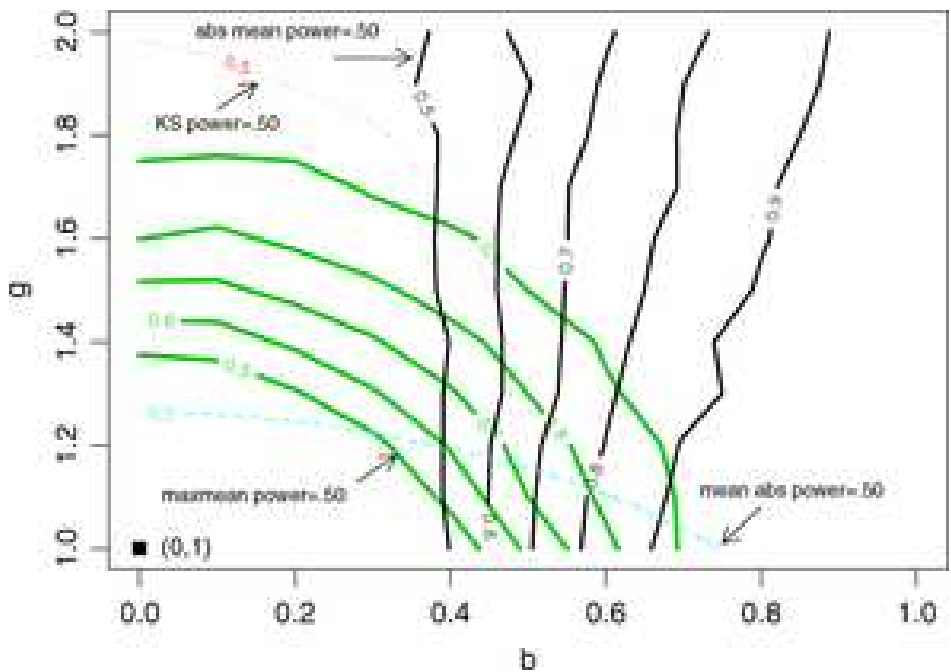

FIG. 5. *Contours of equal power for absolute mean, enrichment statistic $|\bar{z}_\mathcal{S}|$ (jagged vertical line) mean of absolute value (curving broken line), and maxmean statistic (heavy curving lines); for testing $H_0 : z_i \stackrel{\text{ind}}{\sim} N(0,1)$ vs $N(b, g^2)$, $i = 1, 2, \ldots, m = 25$, at level 0.95; "ks" is 0.5 power contour for the Kolmogorov–Smirnov test, dominated by maxmean test with 0.5 contour everywhere closer to null point $(0, 1)$. Using mean absolute value, average over $\mathcal{S}$ of $|z|$ has very little power versus shift alternatives, curved dashed line.*

Using a similar setup, the left panel of Figure 6 compares the power of the maxmean statistic to that of the individual gene scores $z_i$. We see that the maxmean statistic has a strong power advantage. In the right panel we have changed the setup slightly, with half of the genes in the gene-set generated as $N(b, g^2)$, and the other half as $N(-b, g^2)$. The increase in power enjoyed by the maxmean statistic is slightly less than in the left panel, but it is still substantial.

We applied the maxmean statistic to the $p53$ data with the catalog of 522 gene-sets in Subramanian et al. (2005). Restandardization, (2.12), noticeably improved results, though the effect was much less dramatic than in Figure 4. Applying the FDR algorithm to the restandardized $p$-value (2.13), cutoff level $q = 0.10$, yielded 8 significant gene-sets, Table 1, in close agreement with Table 2 of Subramanian et al. (2005).

Significance is defined in a two-sided testing sense here, but we can use $\bar{s}_{\mathcal{S}}^{(+)}$ and $\bar{s}_{\mathcal{S}}^{(-)}$ in (3.11) separately. Gene-sets 4 and 7 in Table 1 have $\bar{s}_{\mathcal{S}}^{(-)}$ predominating, indicating positive enrichment in the mutant class, while the other six gene-sets had bigger $z$-values overall in the normal class.

NOTE. Restandardization formula (2.12) was applied *separately* to $\bar{s}_{\mathcal{S}}^{(+)}$ and $\bar{s}_{\mathcal{S}}^{(-)}$ in (3.11), and then combined to give

$$(3.13) \qquad S_{\max}^{**} = \max\{\bar{s}_{\mathcal{S}}^{(+)**}, \bar{s}_{\mathcal{S}}^{(-)**}\}.$$

Figure 7 displays the top eight pathways, showing the individual gene scores $z_i$ in each pathway. The inset at the top shows the estimated FDR

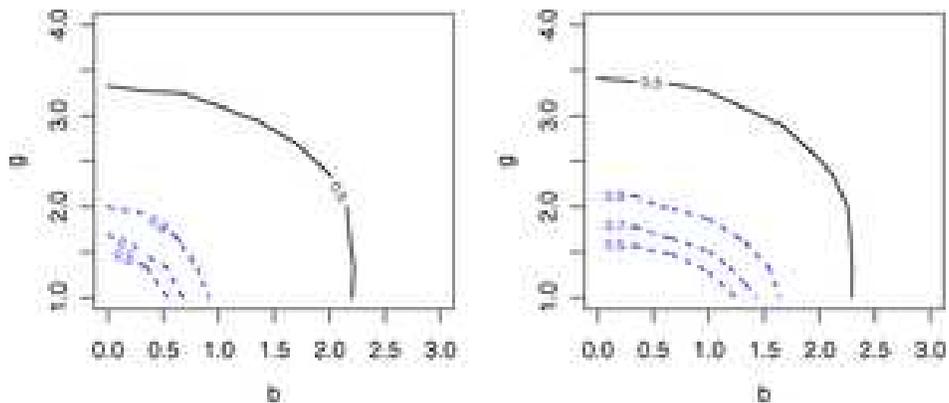

FIG. 6. *Contours of power 0.5 (solid lines) for individual gene scores $z_i$ and power $0.5, 0.7 \ldots 0.9$ (broken lines) for the maxmean statistic; for testing $H_0 : z_i \overset{\text{ind}}{\sim} N(0, 1)$ vs $N(\pm b, g^2), i = 1, 2, \ldots, m = 25$. In the left panel, all genes in the gene-set are shifted by $b$; in the right panel, half of the genes are shifted by $b$ and half by $-b$.*

Table 1
*Significant gene-sets, FDR cutoff 0.10, for p53 data, Subramanian et al. (2005); maxmean statistic restandardized*

| | | | |
|---|---|---|---|
| 1. | $p53$ hypoxia* | 5. | $p53$ Up* |
| 2. | $p53$ pathway* | 6. | hsp27 pathway* |
| 3. | radiation sensitivity* | 7. | n.g.f. pathway |
| 4. | ras pathway* | 8. | $SA$ $G1$ and $S$ phases |

*Indicate significant gene-sets listed in Subramanian et al. (2005), Table 2. They list the n.g.f. pathway as almost significant.

for an individual gene analysis. We see that the cutoff for an FDR of 0.10 is about $-4.7$ for negative scores, while the FDR for positive scores never gets down to 0.10. Hence, an analysis of individual genes would report as significant the top few genes in the negative gene-sets, and none of the genes in the positive gene-sets. Note that location shift alternatives look quite reasonable here.

The values of $(\text{mean}_s, \text{stdev}_s)$ and $(\text{mean}^*, \text{stdev}^*)$ used in (2.12) were based on the catalog of 15,059 genes listed (with multiple occurrences) in the 522 gene-sets, rather than the 10,100 genes in the original $p53$ data set. Table 2 shows considerable differences between the two choices for $(\text{mean}_s, \text{stdev}_s)$, though not for $(\text{mean}^*, \text{stdev}^*)$, with $s = s^{(+)}(z)$, the comparison for $s^-(z)$ being similar. When a catalog of gene-sets $\mathcal{S}$ is simultaneously under consideration for enrichment, the catalog choice of $(\text{mean}_s, \text{stdev}_s)$ is preferable from the row randomization viewpoint.

The Poisson Selection Model (3.1) and (3.2) follows the Randomization Model of Section 2 in that the set $\mathbf{z}$ of $N$ original $z$-values is treated as a fixed ancillary, with only the selection process for $\mathcal{S}$ considered random. We could, instead, think of $z$-values selected for inclusion in $\mathcal{S}$ as coming from a different underlying distribution, for example,

$$f_\beta(z) = e^{\beta' s(z) - \log T_\beta} f_0(z), \tag{3.14}$$

where $f_0(z)$ is the baseline density for nonenriched $z$-values, $\beta$ and $s(z)$ are as in (3.1), and $T_\beta = \int_{-\infty}^{\infty} e^{\beta' s(z)} f_0(z)\, dz$. If the members of $\mathcal{S}$ are selected independently according to (3.14), then the most powerful test for null hypothesis $\beta = 0$ is based on $\bar{s}_\mathcal{S} = \sum_\mathcal{S} s(\mathbf{z})/m$, just as in (3.6).

This point of view is more congenial to the Permutation Model of Section 2, but there is really not much practical difference from (3.1). Let $\widehat{F}_0$ denote the empirical distribution

$$\widehat{F}_0[a,b] = \#\{z_i \text{ in } [a,b]\}/B, \tag{3.15}$$

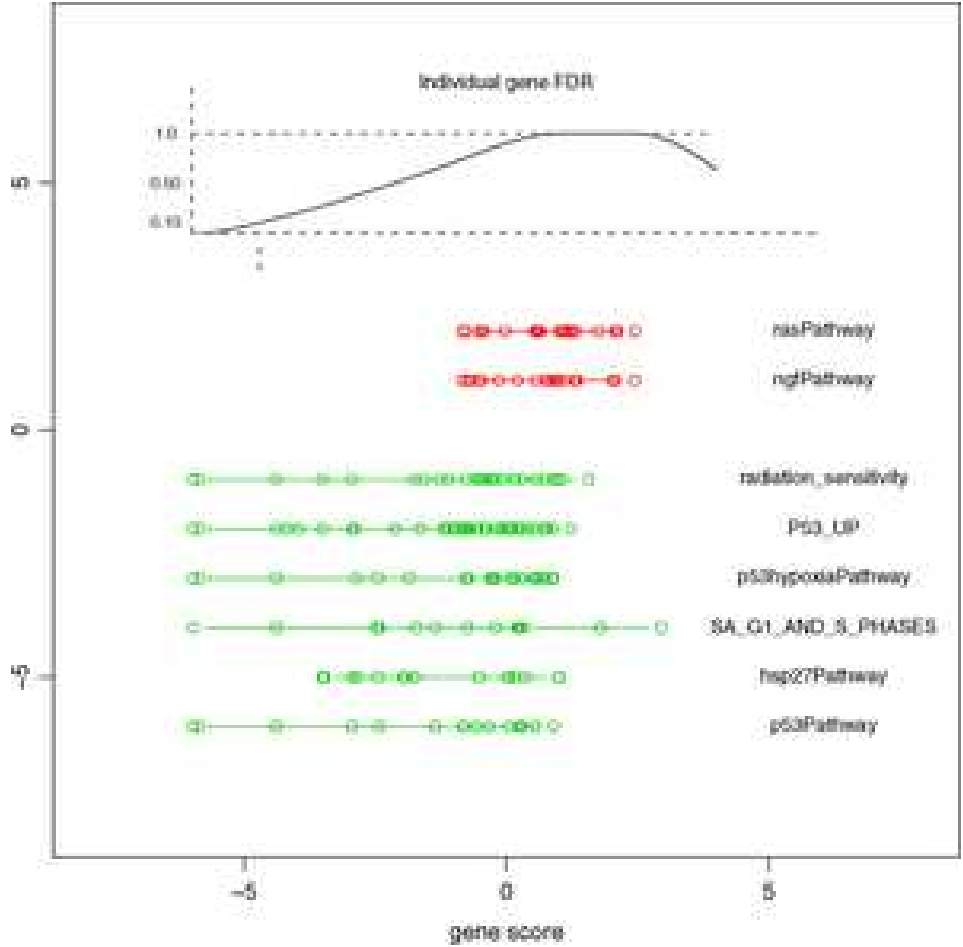

FIG. 7. *p53 data: the top eight pathways showing the individual gene scores $z_i$ in each pathway. The inset at the top shows the estimated FDR for an individual gene analysis.*

TABLE 2
*Comparison of (mean, stdev) values for $s = s^{(+)}(z)$, (3.10); for original list of 10,100 genes in p53 data, versus 15,059 genes listed with multiple occurrences, in catalog of 522 gene-sets. Latter is preferred choice for restandardization formula (2.12)*

|  | **mean$_s$** | **(stdev$_s$)** | **mean$^*$** | **(stdev$^*$)** |
|---|---|---|---|---|
| Original 10,100 genes | 0.451 | (0.628) | 0.406 | (0.589) |
| Catalog 15,059 genes | 0.456 | (0.671) | 0.404 | (0.587) |

and $\widehat{F}_\beta$ the tilted distribution

$$(3.16) \qquad \widehat{F}_\beta[a,b] = \sum \{\mu_i \text{ for } z_i \text{ in } [a,b]\} \Big/ \sum_1^N \mu_i.$$

Then (3.1) and (3.2) give

$$(3.17) \qquad \widehat{F}_\beta[a,b] = \int_a^b e^{\beta' s(z)} \, d\widehat{F}_0(z),$$

just as in (3.14).

More elaborate probability models can be constructed to make, say, the Kolmogorov–Smirnov distance the optimal test statistic. The strongpoint of the maxmean enrichment statistic is its good performance against simple location and scale alternatives, as seen in Figure 5.

**4. Computational issues and software.** The developments in the previous two sections lead to our *Gene Set Analysis procedure*, which we summarize here:

1. Compute a summary statistic $z_i$ for each gene, for example, the two-sample $t$-statistic for two-class data. Let $\mathbf{z}_\mathcal{S}$ be the vector of $z_i$ values for genes in a gene-set $\mathcal{S}$.
2. For each gene-set $\mathcal{S}$, choose a summary statistic $S = s(\mathbf{z})$: choices include the average of $z_i$ or $|z_i|$ for genes in $\mathcal{S}$, the GSEA statistic or, our recommended choice, the *maxmean* statistic (3.11).
3. Standardize $S$ by its randomization mean and standard deviation as in (2.14): $S' = (S - \text{mean}_s)/\text{stdev}_s$. For summary statistics such as the mean, mean absolute value or maxmean (but not GSEA), this can be computed from the genewise means and standard deviations, without having to draw random sets of genes. Note formula (3.13) for the maxmean statistic.
4. Compute permutations of the outcome values (e.g., the class labels in the two-class case) and recompute $S'$ on each permuted dataset, yielding permutation values $S'^{*1}, S'^{*2}, \ldots, S'^{*B}$.

We use these permutation values to estimate $p$-values for each gene-set score $S'$ and false discovery rates applied to these $p$-values for the collection of gene-set scores.

It would be convenient computationally to pool the permutation values for different gene-sets: this is often done in the analysis of individual genes [e.g., in the SAM package Tusher et al. (2001)] and also in the GSEA software. Such pooling reduces the number of permutations needed to obtain sufficient accurate results. However, we have found that in gene-set analysis the permutation distribution for each gene-set depends heavily on both the

ON TESTING THE SIGNIFICANCE OF SETS OF GENES 19

number of genes in the gene-set and their average pairwise correlation (neither of which applies to individual-gene analysis). Hence, we do not pool the permutation values, and so we must carry out at least 1000 permutations to get accurate results. However, with a careful implementation this need not be an obstacle for practical use.

We have written an R language package GSA ("gene-set analysis") for carrying out these computations. This will be freely available along with collections of gene-sets for use with the package.

**5. Simulation comparison of different gene-set summaries.** In this study we simulated 1000 genes and 50 samples in each of 2 classes, control and treatment. The genes were assigned to with 50 gene-sets, each with 20 genes. All measurements were generated as $N(0,1)$ before the treatment effect was added. There were five different scenarios:

1. all 20 genes of gene-set 1 are 0.2 units higher in class 2,
2. the 1st 15 genes of gene-set 1 are 0.3 units higher in class 2,
3. the 1st 10 genes of gene-set 1 are 0.4 units higher in class 2,
4. the 1st 5 genes are 0.6 units higher in class 2,
5. the 1st 10 genes of gene-set 1 are 0.4 units higher in class 2, and 2nd 10 genes of gene-set 1 are 0.4 units lower in class 2.

In every one of these scenarios only the first gene-set is of potential interest. For each scenario, we carried out 200 permutations and report the estimated tail probability $\mathrm{Prob}(S'^* > S'_1)$, with small values being good. Here $S'$ is the restandardized version of a summary statistic, and the quantity $\mathrm{Prob}(S'^* > S'_1)$ is the observed $p$-value for the first gene-set.

We compared five different summary statistics: the mean of $z_i$, the mean of $|z_i|$, maxmean, GSEA and GSEA applied to $|z_i|$. Everything was repeated 20 times, and the mean and standard error of the tail probability over the 20 simulations are reported in Table 3. While the maxmean can be beaten by a small margin by the mean or the mean of the absolute values in the one-sided or two-sided scenarios, respectively, it is the only method with consistently low $p$-values in all five scenarios.

**6. Examples.**

6.1. *Survival with kidney cancer.* In this example we apply the gene-set methodology to a dataset with censored survival outcome. Zhao et al. (2005) collected gene expression data on 14,814 genes from 177 kidney patients. Survival times were also measured for each patient. The data were split into 88 samples to form the training set and the remaining 89 formed the test set.

TABLE 3
*Results of simulation study: p-values for the first gene-set, for five different summary statistics (columns) and five different scenarios (rows, described in the text). The first four scenarios are one-sided shifts, while the last one is two-sided. We see that the maxmean statistic is the only one with consistently low p-values in all five scenarios*

|      | mean  | mean.abs | maxmean | GSEA  | GSEA.abs |
|------|-------|----------|---------|-------|----------|
| (1)  |       |          |         |       |          |
| mean | **0.014** | **0.133** | **0.012** | **0.032** | **0.192** |
| se   | 0.008 | 0.038    | 0.005   | 0.017 | 0.060    |
| (2)  |       |          |         |       |          |
| mean | **0.005** | **0.035** | **0.002** | **0.016** | **0.074** |
| se   | 0.003 | 0.014    | 0.001   | 0.008 | 0.034    |
| (3)  |       |          |         |       |          |
| mean | **0.014** | **0.032** | **0.002** | **0.031** | **0.057** |
| se   | 0.008 | 0.015    | 0.001   | 0.018 | 0.032    |
| (4)  |       |          |         |       |          |
| mean | **0.074** | **0.081** | **0.014** | **0.069** | **0.037** |
| se   | 0.032 | 0.050    | 0.007   | 0.038 | 0.014    |
| (5)  |       |          |         |       |          |
| mean | **0.587** | **5e-04** | **0.018** | **0.233** | **0.011** |
| se   | 0.106 | 4e-04    | 0.008   | 0.063 | 0.009    |

We computed the Cox partial likelihood score statistic $z_i$ for each gene, and used this as the basis for gene-set analysis using the maxmean statistics $S_k$. Since there are separate training and test sets, we examined the concordance between the scores $S_k$ for a given gene-set between the training and test sets. Denote the training and test set scores by $S_k^{tr}$ and $S_k^{te}$. We measured the concordance by setting a cutpoint $c$ and defining a test set value with $|S_k^{te}| > c$ to be a "true positive"; otherwise it is called a "true negative." Then we applied the same cutpoint to the training set scores $S_k^{tr}$ and then computed the false positive and false negative rates. We applied this same idea both to scores for the individual genes and also to gene-set enrichment analysis, based on the Cox scores. The cutpoint $c$ was varied to generate the curves of false positive and false negative rates in Figure 8. We see that both the maxmean statistic and GSEA statistics have considerably lower false positive rates, for a given false negative rate, than the individual gene analysis. The maxmean statistic also enjoys an advantage over GSEA, which is not surprising given the results from our power studies. In this case it seems that a list of significant gene-sets based on the maxmean statistic will tend to be more reproducible than a list of significant genes from an individual gene analysis.

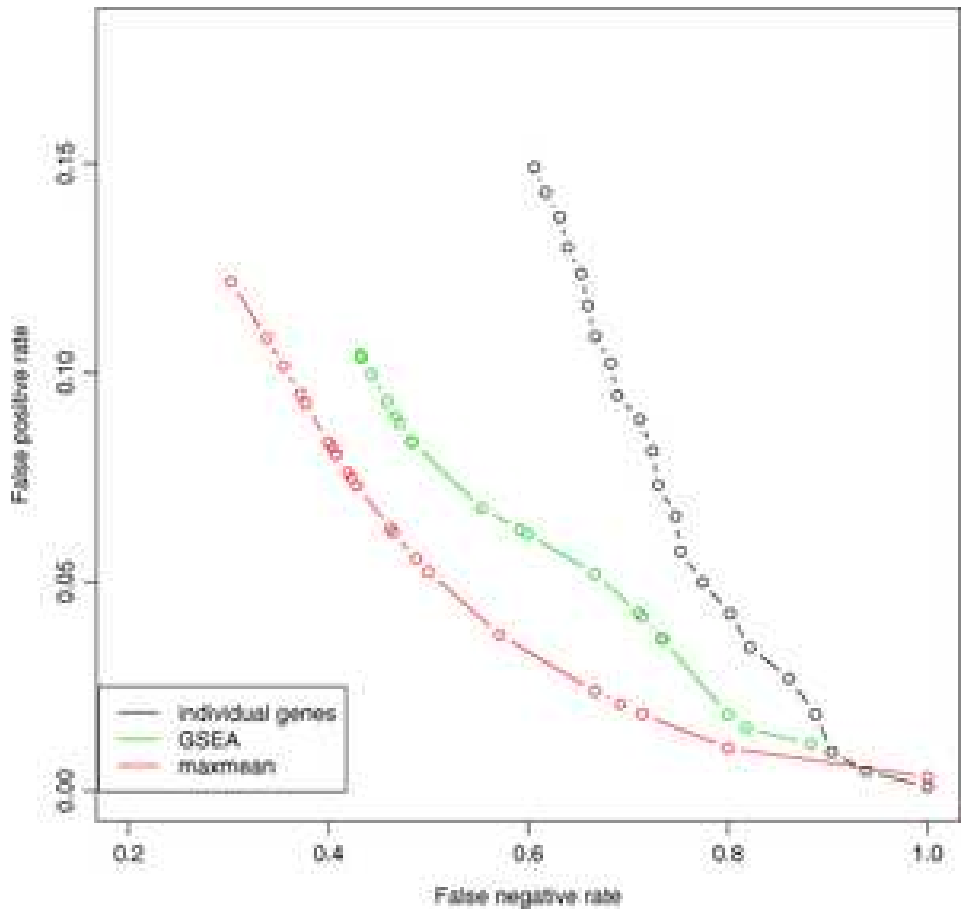

Fig. 8. *Kidney cancer microarray data: false positive and false negative rates for three methods, using the test set values for each method to define the "true positive" and "true negative" for genes or gene-sets.*

6.2. *Comparisons of different gene-set collections.* The gene-sets used in the examples of this paper are the collection of 522 gene pathways developed by Subramanian et al. (2005). This is the C2 "functional" collection in their *Molecular Signature Database*. The same authors provide three additional collections: C1 (chromosomal location, 24 sets), C3 (motif-based, 319 sets ) and C4 (chromosomal location, 427 sets). The Stanford Microarray Database has at least two more collections in their "synethetic gene" database that are potentially useful for this purpose: tissue types (80 sets) and cellular processes (22 sets). We applied our gene-set analysis technique using the *maxmean statistic* to both the p53 and kidney cancer datasets, using each of these six collections of gene-sets. The results are shown in Figure 9. In both cases the FDR results in the left tail are shown for gene-sets with negative maxmean statistics. The FDRs for the right tail were high for both datasets. We see that different collections exhibit the lowest FDRs in the two studies: *motif*, *chromosomes* and *pathways* for the p53 data and *correlation* for the kidney cancer data. This suggests that, for a given expression dataset, systematic "mining" of different gene-set collections could be a useful exercise for discovering potentially interested gene-sets for further investigation.

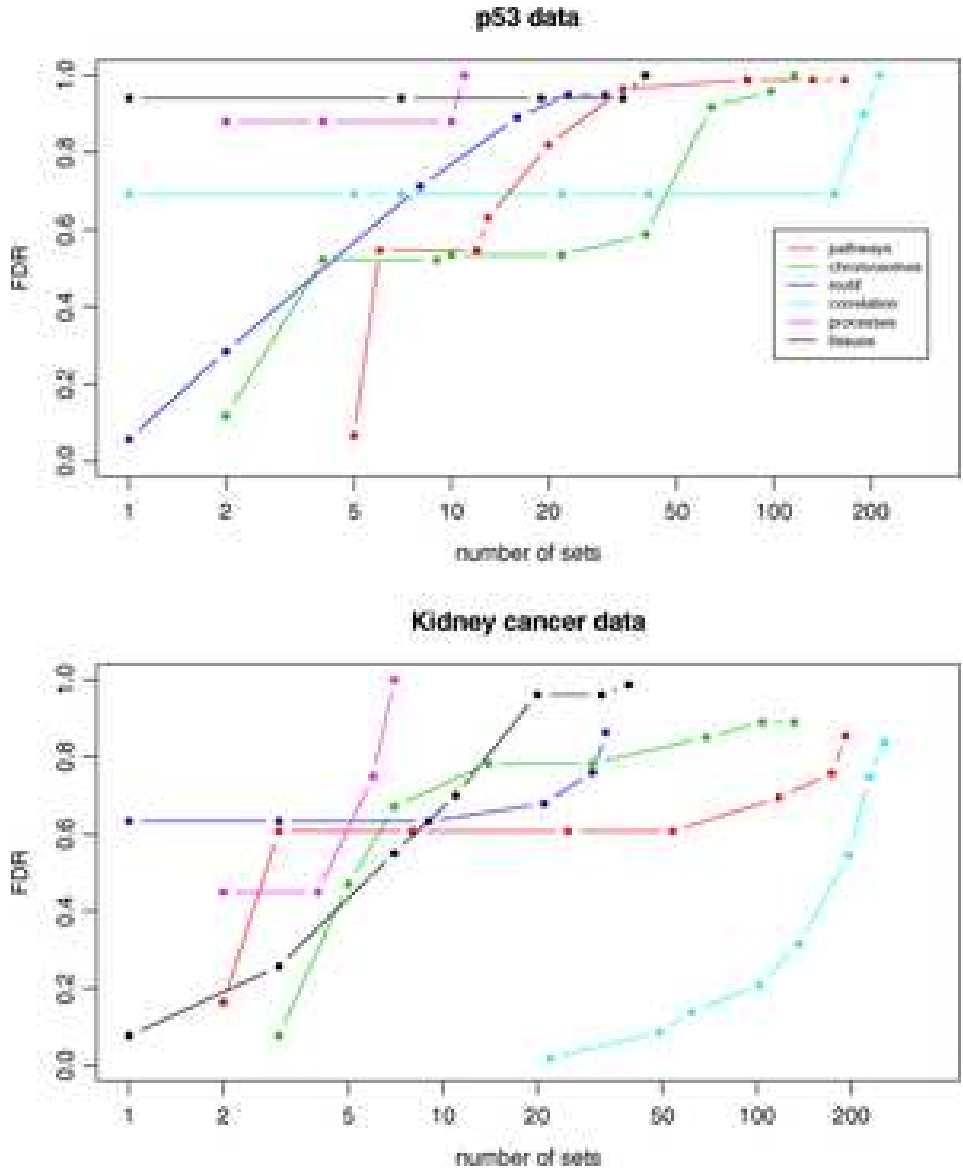

Fig. 9. *Gene set analysis applied to the p53 mutation data and kidney cancer data, comparing the estimated FDRs over six different collections of gene-sets. In both cases the FDR results in the left tail (negative gene-set scores) are shown. The FDRs for the right tail (positive gene-set scores) were high for both datasets.*

**Acknowledgments.** We would like to thank the authors of GSEA for their assistance, especially Pablo Tamayo for considerable efforts in debugging our implementation of the GSEA procedure. We would also like to thank Trevor Hastie and Rafael Irizarry for helpful comments and suggestions. Finally, we would like to thank the Associate Editor and two referees for helpful comments that led to improvements in this paper.

Departments of Statistics,
and Health, Research and Policy
Sequoia Hall
Stanford University
Stanford, California 94305-4065
USA
E-mail: brad@stat.stanford.edu

Departments of Health,
Research and Policy, and Statistics
Stanford University
Stanford, California 94305-5405
USA
E-mail: tibs@stanford.edu